# Measurement of statistical evidence on an absolute scale following thermodynamic principles

V. J.  Vieland[1,2,3], J. Das[1,2,4], S. E. Hodge[1,2], S-C. Seok[1]

[1]Battelle Center for Mathematical Medicine, The Research Institute at Nationwide Children's Hospital, Columbus OH

[2]Department of Pediatrics, The Ohio State University College of Medicine, Columbus OH

[3]Department of Statistics, The Ohio State University, Columbus OH

[4]Department of Physics and the Biophysics Graduate Program, The Ohio State University, Columbus OH

Corresponding Author:

Veronica J. Vieland, Ph.D.

Battelle Chair in Quantitative and Computational Biology

Director, Battelle Center for Mathematical Medicine

Professor of Pediatrics and Statistics

The Research Institute at Nationwide Children's Hospital

& The Ohio State University

575 Children's Crossroad

Columbus OH 43215

Phone: (614) 355-5651

Fax: (614) 355-6671

Email: Veronica.Vieland@NationwideChildrens.org

**ABSTRACT** Statistical analysis is used throughout biomedical research and elsewhere to assess strength of evidence. We have previously argued that typical outcome statistics (including p-values and maximum likelihood ratios) have poor measure-theoretic properties: they can erroneously indicate decreasing evidence as data supporting an hypothesis accumulate; and they are not amenable to calibration, necessary for meaningful comparison of evidence across different study designs, data types, and levels of analysis. We have also previously proposed that thermodynamic theory, which allowed for the first time derivation of an absolute measurement scale for temperature (T), could be used to derive an absolute scale for evidence (E). Here we present a novel thermodynamically-based framework in which measurement of E on an absolute scale, for which "one degree" always means the same thing, becomes possible for the first time. The new framework invites us to think about statistical analyses in terms of the flow of (evidential) information, placing this work in the context of a growing literature on connections among physics, information theory, and statistics.

## INTRODUCTION

In previous work we have argued for the importance of rigorous methods for measuring statistical evidence in biomedical research and other arenas, and we have proposed using the mathematical foundations of thermodynamics as a template for achieving this end (Vieland 2006; Vieland 2011; Vieland and Hodge 2011). Here we complete the first stage of development of a thermodynamically based methodological framework in which evidence becomes measurable on an absolute scale for the first time.

By an absolute scale we mean one for which the unit of measurement *means the same thing* across the range of the measurement scale (just as an increase of one foot of length means the same thing, up to given measurement precision, regardless of whether that one foot is added to one foot or to one hundred feet), and across disparate applications (just as one foot means the same thing whether measuring lengths of rope or lumber). It is worth pointing out from the outset that this understanding of "absolute" goes beyond requiring an absolute 0 point for the measurement scale, which is sometimes misunderstood as the basis on which the Kelvin thermometric scale is said to be absolute. (The term "absolute" is also used differently in representational measurement theory.) The Kelvin scale is clearly and by design absolute in the same sense in which we are using the term (Chang 2004).

Perhaps the most intellectually challenging aspect of our measurement problem is simply articulating what is meant by *means the same thing* for evidence measures. For instance, many statisticians use the p-value to indicate the degree of evidence against a null hypothesis. But how would we go about establishing whether a change of, say, 0.001 *means the same thing* relative to a baseline p-value of 0.01 and relative to another baseline p-value of 0.0001? Other statisticians use the log likelihood ratio (LR) as a measure of evidence, but the same conundrum arises. Not only is it unclear whether a change in the log LR from 1 to 2 *means the same thing* as a change from 10 to 11, or perhaps from 10 to 20; it is moreover not even clear how we would go about investigating the matter. We lack a methodological framework in which asking for the meaning of *means the same thing* is a well-posed question. But surely measurement of any quantity, including evidence, requires at a minimum that the unit of measurement always does in fact mean the same thing.

It is noteworthy that the abstract and obscure nature of this measurement quandary is directly parallel to that faced and overcome by Kelvin himself in the mid-19th century, a stroke of historical luck that we gladly exploit for our own purposes. We develop a methodological framework in the context of which it becomes possible to articulate a definition of *means the same thing* in measuring statistical evidence, based on what is meant by this phrase in reference to degrees of temperature in thermodynamics. In essence our approach is extremely simple: we reinterpret fundamental thermodynamic quantities in terms of evidential rather than physical processes, then port the mathematical underpinnings of thermodynamics over to the evidential side *en bloc*. The quantity denoted by T in thermodynamics – temperature measured on the absolute scale – is in this framework translated as the evidence E, which is then intrinsically measured on the absolute scale, in exactly the same sense in which T is.

Various connections have of course been previously made between thermodynamics and certain aspects of statistical inference; see, e.g., (Shannon 1948; Jaynes 1957; Cox 1961; Kullback 1968; Soofi 2000; Caticha 2003), primarily based on the power of one particular concept – entropy – in multiple settings. We view our own results as extensions of Jaynes' research into connections among Shannon entropy, statistical mechanical entropy, thermodynamic entropy, and statistical inference. We wonder, however, why Jaynes himself never asked the following question: If crucial underlying quantities such as entropy tie statistical inference and thermodynamics together, and given the central role of temperature in thermodynamics, what then is the analogue of temperature on the inferential side? Our answer represents, as far as we know, the first proposal to interpret statistical evidence as the direct mathematical analogue of temperature, and moreover, of temperature as measured on the absolute scale.

We focus in this paper on one very simple statistical model. Suppose we are interested in the evidence that a certain coin is either fair or biased towards landing tails. We perform one or more experiments, each time tossing the coin *n* times and observing the number of times $x$ that the coin lands heads (where both $n$ and $x$ can vary across experiments). For a given experiment, the data ($n$, $x$) carry information about the coin, which we capture by graphing the likelihood ratio (LR) as a function of P(heads) = θ, fixing θ = ½ in the denominator of the LR. Data from subsequent experiments change the LR graph, for instance, causing the height of the graph at its maximum to increase. We consider the evidence itself as a property of the graph, and seek to define this property in a way that is amenable to absolute measurement, while preserving our underlying feel for what we mean by *evidence.*

Our starting point is what may at first appear as a relatively obscure analogy with the elementary heat engine of thermodynamics, more specifically, with the abstract model of the workings of a heat engine embodied in the Carnot cycle. The purpose of a heat engine is to convert heat into work, say, using the heat to raise a piston, and similarly we might describe the purpose of the coin-tossing experiment to be the conversion of the information conveyed by new data into movement of the LR graph. But more importantly for our purposes, the Carnot cycle constituted the basis for rigorous understanding of the fundamental dynamics of thermal systems. We therefore consider what is to be gained by considering the statistical experiment as the analogue of the heat engine, in order to rigorously characterize the relationships among evidential analogues of heat, work, and temperature. Thus rather than begin by assuming any particular definition of E, instead, we adopt as our sole methodological constraint that our evidential model must adhere to the equations required to construct and analyze Carnot cycles. The resulting evidential dynamical framework will be justified insofar as it proves fruitful, even if it should remain difficult to relate the idealized workings of an engine to how we ordinarily think about statistical data analyses. Indeed, if this approach does prove to be useful, we may question whether it is based on mere analogy, or possibly, the far deeper connection between thermodynamics and inference championed by Jaynes and many others.

## METHODS AND RESULTS

We develop the central argument of this paper in three steps, as follows. We begin by postulating that a statistical model can be described by an equation (called an equation of state, or EqS), expressed in the form of the EqS for an ideal gas, and that this EqS defines a system that is governed by analogues of the 1$^{st}$ and 2$^{nd}$ Laws of Thermodynamics. In section (1) below, we reparameterize a simple binomial model to illustrate. *This establishes that it is <u>possible</u> to represent a statistical system in thermodynamic terms.* We then (2) demonstrate that the quantity corresponding to temperature in this system, which we call E (for evidence), exhibits behavior that coincides in multiple respects with our intuitive feel for the behavior of statistical evidence. *This establishes that it is <u>plausible</u> to use such a system to define statistical evidence.* Finally, (3) we explicitly address the meaning of *means the same thing* for the unit of measurement for E. We show how this formalism leads us to think about statistical systems in fundamentally new ways, and *we posit that doing so will prove <u>useful</u>*, just as early thermodynamic investigations led to new and productive ways of thinking about heat. Throughout we attempt to motivate the analogy between specific thermodynamic quantities and their evidential counterparts heuristically; in Appendix A we present an alternative derivation of fundamentals of the new framework, which is mathematically more direct but also more abstract.

**(1) *Possibility of the Formalism*** We start with a simple binomial coin tossing model, with *x* heads out of *n* tosses, and probability θ that the coin lands heads. We consider the hypothesis contrast "coin is biased" versus "coin is not biased" (θ=0.5).[1] We stipulate that the LR (Eq. 1, Table 1), captures all the evidential content of the system, insofar as the addition of anything beyond the LR would not change the evidence conveyed by given data with respect to the stated hypotheses. This stipulation is fundamental to our framework; for related discussion see (Hacking 1965; Edwards 1992; Royall 1997). While we generally think of changes in the LR graph as resulting from changes in the data ($n$, $x$), we can equally well consider the permissible transformations of the graph *per se*, where what is permissible is any state of the system that conforms to the underlying binomial LR equation as the data change. In section (2) below, we treat both $n$ and $x$ as continuous (non-negative) numbers, for reasons that will become clear in section (3).

We begin with the postulate that the EqS for this system can be described by the equation describing an ideal gas (Eq. 2, Table 1). Our task, therefore, is to see if we can find a reparameterization of the statistical system that enables us to express it in the form PV=RT. This involves deriving the analogues of the constituent quantities and demonstrating internal consistency for the resulting EqS in statistical terms. (Here we omit the usual factor N = number of moles of matter, or equivalently, we use as our template the EqS for a single mole of matter; see the Discussion for further comments.)

As the evidential analogue of volume, $V_E$, we simply use the area under the LR curve (Eq. 3, Table 1). (In a two-parameter model, this would be an actual volume.) Here we consider a one-sided case, in which the quantity in Eq. 3 is integrated from θ=0 to 0.5 (with $x \leq n/2$). While this complicates the arithmetic somewhat, it will simplify some aspects of the presentation. The one-sided case also plays a special role in statistical genetics, and we have investigated its properties in detail elsewhere (Hodge and Vieland 2010). We have established previously that at least in statistical genetic applications, $V_E$ exhibits "thermoscopic" behavior, in the sense that it properly goes up (or down) as the evidence goes up (or down), and in a manner unmatched by other standard approaches to quantifying evidence, such as the p-value or the maximum LR (Huang and Vieland 2001; Bartlett and Vieland 2006; Huang and Vieland 2010). Thus we view $V_E$ as a direct analogue of physical volume V in elementary thermodynamic systems, in that like V, $V_E$ can serve as a thermoscopic indicator of increasing or decreasing temperature.[2]

---

[1] The form of this contrast is fundamental to the formalism and relates to the role of information entropy in the system; see below.

In thermodynamics, if we impose the 1[st] Law (conservation of energy) and 2[nd] Law (impossibility of a perpetual motion machine) on an ideal gas, then we can also express the EqS in terms of the entropy S (Fermi 1956 (orig. 1936)) (Eq. 4 in our Table 1), where S is defined up to an (additive) constant. In physics, $C_V$ can take on one of three values: 1.5R, 2.5R, 3R, for monatomic, diatomic, and polyatomic gasses respectively. Here we fix $C_V$=1.5R throughout. See the Discussion and Appendix A for further comments on the evidential analogues of the 1[st] and 2[nd] Laws; see also Palacios (quoted in (Krantz, Luce et al. 1971, p. 456)) for discussion of the thermodynamic constants.

From Jaynes (Jaynes 1957) and Kullback (Kullback 1968) we have a direct connection between the thermodynamic quantity S and the Shannon information, and between the Shannon information and statistical information (Fisher information, etc.). We therefore deploy a definition of S based on Shannon. This will be the "hook" that enables us to connect the binomial LR to a thermodynamic representation of the system. Specifically, we define the evidential analogue of S, which we denote $S_E$, as a particular form of relative entropy, viz., the information entropy evaluated at the constrained maximum entropy (MaxEnt) value of the binomial probability θ, which is in this case simply the maximum likelihood estimator $x/n$, relative to the information entropy evaluated at the unconstrained MaxEnt value θ = ½ (Eq. 5, Table 1). As in physics, $S_E$ is defined only up to an additive constant. See also (Krantz, Luce et al. 1971) for an axiomatic basis for entropy as a mathematical rather than physical quantity. We note that the property of reversibility is a consequence of this definition: $S_E$ depends only on the current state of the system, regardless of the path taken on a $P_E V_E$ diagram (see below) to arrive at that state from some other state; and because we have defined it in terms of MaxEnt states, assuming that MaxEnt also corresponds to statistical equilibrium, the evidential system developed here can be said to be in an equilibrium state at all times.

---

[2] In previous work we have referred to a quantity closely related to $V_E$ as the Bayes Ratio (BR, see e.g. Vieland et al. (2011)). The BR differs from $V_E$ insofar as it multiplies the LR by a prior distribution, generally taken to be uniform. Note that $V_E$ also differs from the Bayes Factor in general contexts, insofar as the integral of $V_E$ is taken over the LR as a whole. In the present case, however, this distinction is moot as the denominator is not a function of $\theta$.

We can now *derive* the remaining variables needed to completely describe the system. Substituting $S_E$ and $V_E$ into Eq. 4 and exponentiating, we obtain an evidential analogue of T, which we call E (Eq. 6, Table 1). Note that E in Eq. 6 is always non-negative. (Because $S_E$ is defined with respect to the MaxEnt state, we expect this definition of E to relate directly to the parameter β of the Boltzmann distribution.) Finally, substituting E as expressed in Eq. 6 into Eq.2 (Table 1), we obtain the evidential analogue of the pressure P, $P_E$ (Eq. 7, Table 1).

Thus we have successfully reparameterized the initial binomial system, expressed in terms of (LR, $n$, $x$), into an explicit EqS describing the LR graph in terms of ($V_E$, $P_E$, E), while tacitly imposing analogues of the 1st and 2nd Laws of Thermodynamics. We can now describe allowable transformations of the system in terms of any change in the LR graph that conforms to the underlying binomial EqS, without explicitly referencing changes in data. This shift from the usual statistical perspective is critical, because it allows us to view changes in the LR graph in terms of the in- or out-flux of *evidential information*, rather than in terms of changes in data (see section 3 below). The only remaining unknown in this system is the constant R. For the present, we set R = 1 (again, see (Krantz, Luce et al. 1971)).

**(2)** ***Plausibility of the Formalism*** In this section we look at features of the system to see how they relate to what we would ordinarily mean by statistical evidence. As just noted, a crucial aspect of this formalism is that we can consider transformations of the system directly in terms of changes in the LR graph, rather than in terms of changes in the data. However, in order to establish the plausibility of the system, it is desirable to be able to appeal to statistical intuition, and this requires considering the behavior of the graph in familiar terms, as a function of $(n, x)$. Thus while we continue to assume that we are working with the EqS for (on the physical side) a single mole of matter, or what is probably best viewed as a fixed quantity of data, we will plot results as a function of changes in $n$. The rationale lies in the transitional nature of this section, illustrating properties of the system in familiar statistical terms, *en route* to thinking of changes as influenced not by the influx of data, but rather, by the influx of evidential information (again, see section 3).

A feature of this system that plays an important role in subsequent discussion is that it captures evidence either in favor of the numerator of the LR or in favor of the denominator.[3] However, there is no particular "mark" designating which is which, that is, no fixed value below which the evidence always favors the denominator. On the other hand, for given $n$, there exists a value of $x$ at which E is minimized (Figure 1(a)). We call this $(n, x)$ pair, which depends on E, the TRansition Point (TrP(E)).[4] To the left of TrP(E) the evidence favors "coin is biased," and to the right of TrP(E) the evidence favors "coin is not biased."

---

[3] We have argued elsewhere that this is a critically important feature of any evidence metric (see e.g. Vieland (2006)).

The behavior shown in Figure 1(a) makes sense from a statistical point of view. Bearing in mind whether one is looking to the left or right of TrP(E), this behavior conforms to what we mean by evidence in terms of two basic properties: (*i*) for fixed *x/n*, as *n* increases E increases; (*ii*) for fixed *n* and reading from left to right on the *x/n* axis, E decreases up to TrP(E), then increases. Figure 1(b) shows *x/n* at TrP(E) over a broader range of *n*. This plot also illustrates sensible behavior: as *n* increases, TrP(E) moves towards $x = n/2$; i.e., for very large *n*, even a small deviation from $x = n/2$ constitutes evidence against $\theta = 0.5$.[5] One thing that may trouble statisticians is the increase in E at TrP(E) as a function of n. We return to this point in the Discussion. Figure 1(c) shows a view of E as a function of $(n, x)$ in three-dimensional space for *n* up to 1000, again showing consistent statistical behavior. We note also that as Figures 1(a) and 1(c) illustrate, for given *n*, it is possible to have far stronger evidence in favor of "coin is biased" than in favor of "coin is not biased." This reflects the fundamental asymmetry in the hypotheses, with the denominator specifying a single value ($\theta=\frac{1}{2}$), and the numerator allowing for any value $\theta < \frac{1}{2}$.

In aggregate, the behavior illustrated in Figure 1 speaks to the *plausibility* of the current model as a representation of a statistical system. Moreover, with the model in hand, we can now begin to think of the system in ways that are unfamiliar to statisticians. For instance, the contours of Figure 1(c) describe what could be called the "isotherms" of the system, that is, the set of ($n$, $x$) values corresponding to a single value of E (Figure 2(a); Figure 2(b) shows these same isotherms plotted against $x/n$). These isotherms also behave in accordance with our intuitions regarding evidence. E.g., consider ($n = 10$, $x = 0$) (not shown on Figure). Whatever the value of the evidence for that point, we know that as $n$ increases while holding $x$ fixed, the evidence must go up. This implies that, were we to hold the evidence fixed and increase $n$, $x$ would have to increase in order to compensate. Conversely, if we hold $n$ fixed and increase $x$, the evidence will diminish, which implies that if we were to hold the evidence constant while increasing $x$, $n$ would have to go up to compensate. We do not normally think of the dynamics of statistical systems in these terms, precisely because outside of this framework, we have no way to hold the evidence constant in our mind's eye while allowing the LR graph to change. Given a way to establish these isotherms, however, the behavior illustrated in Figure 2(b) clearly makes statistical sense. Note also that the maximum value of $n$ for each isotherm occurs at TrP(E). Figure 2(c) shows one isotherm plotted as a function of $V_E$ and $P_E$, illustrating consistency between the behavior of E and a physical system.

---

[4] TrP(E) is obtained by solving the following equation for *x*:  $\dfrac{\int_{x}^{\frac{n}{2}} \theta^{x}(1-\theta)^{n-x} \log\left(\dfrac{\theta}{1-\theta}\right) d\theta}{\int_{x}^{\frac{n}{2}} \theta^{x}(1-\theta)^{n-x} d\theta} = \log\dfrac{x}{n-x}$ .

[5] If we were considering values of *x > n/2*, we would also need to distinguish a corresponding set of mirrored sectors on that side of *n/2*. This is one reason it is simpler to restrict attention to the one-sided case, in which we need only consider two sectors of the graphs.

Finally, we consider another quantity unfamiliar in traditional statistical frameworks, namely, the evidential equivalent of physical (mechanical) work. We are interested in the effects of in- or out-flux of evidential information on the LR graph, characterized in terms of the fundamental features of the graph. We posit, therefore, that the quantity $P_E dV_E$ correctly characterizes a transformation of the system from an information-transfer point of view, just as PdV does in physics from an energy-transfer point of view (Fermi 1956 (orig. 1936)). Thus the amount of evidential work done during a given transformation of the system from state A to state B is the quantity $W_E$ (Eq. 8, Table 1). While we are not aware of a familiar statistical analogue of this notion of evidential work, nevertheless, we think that work defined in this way has some intuitively appealing features. In particular, a given change in $V_E$ is more difficult to effect (requires greater information influx, see below) for higher E. As in physics, we distinguish work being done by the system ($W_E$) from work done to the system ($-W_E$). The central importance of this notion of evidential work will become clear in the following section. See also Appendix A, which makes clear that nothing in the underlying mathematics requires an interpretation in terms of *mechanical* work per se.

Note that, because our evidential system follows the ideal gas EqS by design, it follows immediately that the system exhibits the characteristic properties of ideal gases, including Boyle's Law (for fixed E, $P_E$ and $V_E$ are inversely proportional), Charles' Law (for fixed $P_E$, $V_E$ is proportional to E), etc.

**(3)** ***Utility of the Formalism*** We have thus far established that it is *possible* and *plausible* to view a statistical system in thermodynamic terms, with the evidential analogue of T, our quantity E, serving as a measure of statistical evidence. It remains to illustrate why we believe this new framework will prove *useful*. Here we show the utility of the framework for establishing an absolute scale for the measurement of E and a context in which we can explicitly give the meaning of *means the same thing* for one degree on this scale.

Following Kelvin, we make use of the Carnot engine, an abstract mathematical device for considering the relationship between physical heat and work.[6] As we now know, although Carnot himself did not, a heat engine works by converting one form of energy to another; specifically, by converting heat input into mechanical work. The Carnot engine can be pictured as a cylinder containing a fixed quantity of gas, with a movable (frictionless) piston constraining the volume, and with two (infinite) heat reservoirs available for contact with the cylinder, one at temperature $t_2$ and the other at $t_1 < t_2$. Here we use lower case $t$ to indicate that the temperature need not be measured on the absolute scale; any thermoscopic measure will suffice, even one for which the size of the unit changes across different parts of the temperature range. Each cycle of the Carnot engine is decomposed into four distinct strokes, A-D, as follows: (A) Heat is absorbed by the system from the warmer heat reservoir while maintaining constant temperature $t_2$, with corresponding increase in V and decrease in P, for an isothermal change. (B) Heat transfer is stopped (the reservoir is removed from contact) and the system is allowed to continue to expand and cool down to temperature $t_1$, through a process known as adiabatic change. (C) The system is compressed isothermally by application of external force (work) at constant $t_1$, with corresponding reduction in V and increase in P, requiring a corresponding transfer of energy out of the system into the cooler heat reservoir. (D) The system is permitted to continue adiabatically, with no further transfer of energy (no further contact with the second reservoir), until it returns to its initial state. During the first two strokes work is done by the system (positive work), while during the second two strokes work is done to the system (negative work). A mechanical heat engine is effective insofar as more work is done by the system during the first two strokes than is required to be done to the system through external means during the third stroke in order to return it to its initial state, and this happens because less work is required to effect the same transformation (or its reverse) at lower temperatures.

---

[6] Throughout this section we draw heavily on Fermi (1956), where the equations for adiabatic transformations used in the evidential Carnot cycles can also be found.

We too can run *evidential Carnot cycles* using our binomial system.[7] Figure 3 shows two numerical examples of such cycles on a classical PV (more specifically, $P_E V_E$) diagram. These cycles are readily seen to correspond to genuine Carnot cycles: in each cycle, work done during the two adiabatic strokes ($W_B$, -$W_D$) is identical in magnitude but opposite in sign; work $W_A$ done in the first (isothermal) stroke is greater in magnitude than the work –$W_C$ done in the third (isothermal) stroke; and $W_C/W_A = e_1/e_2$, where $e$ is the analog of $t$ in the physical engine, that is, any "thermoscopic" measure of evidence. We now posit that evidential work is a transformation of evidential information, just as mechanical work is a transformation of heat in a physical heat engine. In physics, heat $Q$ represents the quantity of energy transferred in or out of the system. Following this convention, we consider *evidential energy*, $Q_E$, to be the amount of evidential information being transferred. In what follows we use $Q_2$ to indicate the evidential energy transferred into the engine during the first stroke at $e_2$, and $Q_1$ to denote the evidential energy transferred out of the engine during the third stroke at $e_1$.

---

[7] The Carnot cycle requires reversibility of all processes, that is, it is an idealized system which never departs from equilibrium. This property does not hold for any actual physical system. Our evidential system appears to be inherently reversible, at least insofar as it is always (instantaneously) in its MaxEnt state.

For a cyclic, reversible transformation, there is no net change in energy from the initial to final state of the system, since by stipulation the system is returned to its initial state. In this case and for an ideal gas, W is a direct measure of heat, and in particular, $W_A = Q_2$. Thus defining the efficiency $\eta$ of the engine as the ratio of the work performed by the cycle to the evidential information absorbed at the higher evidence level (and temporarily omitting subscript "E"s for notational convenience), we have $\eta = (W_A+W_B-W_C-W_D)/W_A = (W_A-W_C)/W_A = 1-W_C/W_A = 1-Q_1/Q_2$. Carnot proved that reversible cyclic engines have maximal efficiency, and by direct application of his reasoning we also have the result that all such evidential engines operating between the evidence levels $e_1$, $e_2$ share the same ratio $Q_1/Q_2$. In other words, the ratio of evidential energies depends only on the ratio of evidence levels, and on nothing else specific to the particular system. Thus we can write $Q_1/Q_2 = f(e_1, e_2)$, where (paraphrasing Fermi) *f is a universal function of the two evidence levels*. From here we can directly follow the arithmetic in Fermi to show that $f(e_1, e_2) = f(e_0, e_1)/ f(e_0, e_2)$, for $e_0$ any arbitrary level of evidence. Fixing $e_0$ at a constant value we can therefore consider $f(e_0, e)$ to be a function of $e$ only, thus, $K f(e_0, e) = g(e)$, where K is an arbitrary constant. This gives, finally, $Q_1/Q_2 = g(e_1)/g(e_2)$.

If we now consider g($e$) to be the evidence E, we see that we have converted our "thermoscopic" measure $e$ onto an absolute scale in the following sense: first, it is immediately apparent that E is on a ratio scale (Hand 2004), hence has an absolute 0; second, while the functional form of g remains to be set, meaning that the size of the degree is not yet specified, a given change in evidential energy will always correspond to the same amount of change in E. This result is entirely independent of the particulars of the engine, and will hold for any reversible cyclic operation.[8] *Figure 3 illustrates the meaning of "means the same thing" in the context of reversible cyclic transformations. As shown, a two-fold change in E corresponds to an identical ratio of evidential work, and thus also an identical ratio of evidential information, across the measurement scale.*

As in physics, we did not define $Q_E$ up front. However, we have now arrived at an understanding of what $Q_E$ must be, namely, the transfer of evidential information in a form that can be converted to $W_E$. Moreover, we immediately have a ready mechanism for measuring the amount of evidential information through the heat-work (energy-work) relationship in an evidential Carnot cycle. In principle then, we can use any "thermoscopic" measure of statistical evidence, such as $V_E$ itself (at least for simple systems such as the one-sided binomial), as an evidence measure, and calibrate it against this absolute scale. Of course, different functions $g$ must be derived for each different type of statistical outcome measure in order to accomplish this; and beyond the context of the highly stylized Carnot cycle, myriad practical issues complicate actual measurement calibration. Additionally, the statistical outcome measure must itself behave thermoscopically, properly tracking up and down with the evidence.

---

[8] In particular, for an ideal gas, $W_C/W_A$ is independent of the heat capacity constant $C_V$ in Eq 4 (Table 1). This is one reason we have not worried about the value of $C_V$ in the preceding development of the theory. See also Hand's (2004) discussion (pages 65-67) of Luce's principles and Krantz et al. (2007) for some relevant measure-theoretic principles.

Thus we have not so much solved the evidence measurement problem, as reformulated it in a way that makes it amenable to solution for the first time. Any remaining skepticism regarding the probability of success in practical applications may be offset by an appreciation of how difficult the corresponding practical problems were in physics (Chang 2004). It is all the more heartening, therefore, that we can now take the ordinary drugstore thermometer wholly for granted.

**DISCUSSION**

We can now restate the central problem we set out to address. We imagine that a given set of data has internal energy, $U_E$. $U_E$ can be thought of as a form of information if that seems more comfortable, but to begin with we did not define $U_E$, nor did we propose any particular device for directly measuring it. (This too parallels the development of thermodynamics (Callen 1985) (p. 461).) We posit only that it is the transfer of evidential energy in or out of a statistical system that causes changes in the LR graph. Having characterized those changes as above, we can now refer to this as the process of an influx (or outflux) of evidential heat $Q_E$ doing evidential work $W_E$. Our task then is to measure that property of $U_E$ that corresponds to what we mean by statistical evidence, which we have done by drawing directly on thermodynamic results to show how to measure E, the evidential analogue of T, on what is demonstrably an absolute scale in the same sense in which the Kelvin scale is an absolute scale for temperature. As in thermodynamics (Chang 2004), we first derived the measurement scale, and only then were we in a position to consider what precisely it is that E is measuring.

Under this new formalism, E is a measure of evidential information, or what we might now preferably call *evidential energy*, with the unit of E having fixed meaning as described above. This opens an entirely new line of research into the manner in which evidential energy is transformed into evidential work as data accumulate, including investigation of relative efficiencies of different mechanisms (statistical models) for so doing. We speculate that this shift in perspective will have numerous ancillary benefits, allowing us to extend Jaynes' program of connecting the laws of inference with the laws governing physical systems through the shared concept of information, or evidential energy. (Note, however, that our framework is silent on some matters of central concern to Jaynes and other Bayesians, such as the proper formulation of priors or rational procedures for rank-ordering beliefs.) First and foremost, however, we envision this paper as merely a prelude to the difficult empirical work of aligning evidence measures across applications in biomedical research, so that results can be properly interpreted across the measurement scale within any one application, across disparate applications, and across different types and levels of analysis.

While it may seem odd, although we have derived an absolute scale for evidence measurement we have in fact not yet specified the units of that scale. In thermodynamics too, the size of the degree is an arbitrary choice, set by convention. Thus we can say that E, as it occurs in all Figures in this paper, is on the absolute scale; yet at the same time, the numerical values of E shown here may not correspond to what we will eventually settle upon as the actual values of the evidence, once a decision has been made regarding the size of the standard measurement unit.[9]

---

[9] One interesting departure from physical systems, at least for simple evidential systems such as the binomial, is the apparent lack of a need for dimensional analysis. We speculate that this relates to the fact that, apart from the evidence itself, the remaining underlying evidential quantities are all themselves immediately on

One important consequence of our results is the necessity to reconsider the role of $n$ in statistical systems. Initially, we had supposed that the accumulation of evidence based on two sets of data, with $n_1 = n_2$ observations respectively, would be the analogue of doubling the number of particles (or moles of matter) in the system. However, the new framework invites us to view $n$ (more precisely, the pair $(n, x)$) as an index of evidential energy, not of the number of particles. Once we see that E measures energy $U_E$ and that in the context of an evidential Carnot cycle this energy can be directly measured through the work $W_E$, it becomes clear that $n$ as it appears in our equations must relate to the units of energy, not matter: $(n, x)$ is what changes as energy flows in and out of the system, and it is changes in $(n, x)$ that correspond to evidential work. This is built into the thermodynamic equations that form the basis of the system. In thermodynamics too, it is possible to effect equivalent changes in macroscopic properties of a system (such as V, P) by either the input of energy or the input of particles (under suitable conditions for each). Indeed, our own previous work on evidential "thermoscopes" has highlighted the distinction between "pooling" data (considering all data as a single data set) and "sequentially updating" the evidence across data sets (Vieland 2006; Hodge and Vieland 2010; Vieland, Huang et al. 2011). In hindsight this appears to relate to the distinction between viewing changes in the system as resulting from changes in $n$ considered as the number of observations ("pooling") or as changes in $n$ considered as an index of energy (sequentially updating results without combining data). Virtually all asymptotic theory in statistical inference is based on pooling as the fundamental operation as n increases, but we have argued elsewhere (Vieland and Hodge 2011) that a cogent measure-theoretic approach to the concatenation of evidence across data sets may require some form of updating. This remains a topic for further research.

---

numerical scales, that is, they are not measurements of anything in the usual sense. Hence their units are implicit and need not be further specified.

We also view it as likely that future work will need to carefully consider the units of $n$ itself. E.g., in human genetics, the unit of observation may be a pair of individuals (e.g., an affected sibling-pair) or an arbitrarily large family, which might correspond to gases of different atomic structures. Another tantalizing connection relates to the characterization of energy in our system by the $(n, x)$ pair, with $x$ being bounded by $n$. Possibly even a simple statistical system would be more appropriately modeled through the Van der Waals EqS (Fermi 1956 (orig. 1936)), which incorporates particle size and particle interaction, with the strength of the latter being constrained by the former. While merely a metaphor at this point, the Van der Waals equation has the further advantage that it allows for  phase transitions, in a way that Eq. 2 itself does not. It could be that the two sides of TrP(E) are best thought of as different "phases," like liquid and gaseous states. This might also relate to the behavior of E as an increasing function of n at TrP(E), which  is a direct consequence of the definition of E (Eq. 6) as a function of $S_E$. Possibly in further elaborations of the model, TrP(E) will prove to be a phase transition boundary point, at which E remains constant regardless of energy input into the system. See also Appendix B, in which we derive a relationship between E and the observed Fisher Information. This relationship is interesting in its own right; moreover, we believe that further consideration may yield insights into the behavior of E at TrP(E).

One further aspect of our system that may strike both statisticians and physicists as implausible is the assumption of evidential analogues of the 1st and 2nd Laws of Thermodynamics, which are generally thought of as essentially physical in nature. In Appendix A we formally define the evidential 1st Law and take preliminary steps towards an evidential 2nd Law. (See(Vieland 2013) for a more complete treatment of this topic.) As in thermodynamics, these laws cannot be proved, but will be vindicated insofar as they allow us to usefully describe and manipulate statistical systems (Van Ness 1983; Callen 1985). Given the coherence of the new evidential framework so far, derived under the methodological assumption of an intimate relationship to the equations of thermodynamics, we believe that further investigation will confirm deep axiomatic connections between statistical systems characterized in terms of evidential information or energy, and the basic laws of thermodynamics.

Theoretical matters aside, the ultimate intent of this project is to be useful in practice. Measurement is not, after all, an end in itself, but rather a *sine qua non* of rigorous scientific activity. When completed, our formalism would not replace other statistical investigations, but it would ideally provide a basis for reporting results of statistical analyses on a unified scale for purposes of meaningful comparison across the range of the measurement scale and across disparate applications. Considerable hard work remains to be done, however, before this framework can be usefully applied to real data analyses. To begin with, we have considered only a simple, discrete distribution. (Discrete distributions may play a special role in evidential systems, just as they do in physical systems.) We will need to establish the validity of corresponding treatments of more complex statistical models involving nuisance parameters, approximating likelihoods, and disparate data structures, along with two-sided LR contrasts, before the framework can be deployed in target areas such as biomedical research. Finally, as a practical matter, we will need ways to indicate whether a given value of E is to the left, directly at, or to the right of TrP(E). This information is readily recovered from the underlying calculations, but forms a kind of meta-data that would need to be associated with the value of E in applications. Again, this might be equivalent to wanting to know whether a substance is in its liquid or gaseous state, which depends not just on temperature, but also on the volume and pressure of the system.

The practical challenges to deploying the new framework are substantial, but no more so than the challenges physicists faced in making practical use of the Kelvin thermal scale itself (Chang 2004). Just as we have exploited the well-established foundations of thermodynamics to develop the new evidential model, we believe that we can look to the experiences of experimental physics to facilitate translation of the theory into practice. Here too we may enjoy the benefits of standing on the shoulders of giants.

## Appendix A: Evidential versions of the 1$^{st}$ and 2$^{nd}$ Laws of Thermodynamics

(Material in this Appendix draws on(Vieland 2013).) Consider the LR graph describing a set of data. In any reasonable statistical treatment, this graph changes only if we alter the data; otherwise we are working with something other than the LR under its usual definition. In particular, the consideration of new data will in general result in a transformation of the LR graph from its initial state (for initial data) to some new state (for the initial and new data considered together). Even prior to considering the nature or mechanism of that transformation, we can insist that it must appropriately reflect the effects of the new data and nothing but the effects of the new data. As a statistical matter, this hardly seems worth stating as some kind of law, but it can be seen to be the analog of an important physical principle:

> ***Conservation Principle*** *(physical): The variation in energy of a system during any transformation is equal to the amount of energy that the system receives from its environment.*

Fermi (1956, p. 11) gives this as an informal statement of the 1$^{st}$ Law of thermodynamics (1$^{st}$ Law).

He then says: "In order to give a precise meaning to this statement, it is necessary to define the phrases "energy of the system" and "energy that the system receives from its environment during a transformation." In this Appendix we articulate the corresponding definitions on the evidential side, in order to precisely state an evidential analogue of the $1^{st}$ Law. We then comment briefly on implications of our work thus far regarding an evidential analogue of the $2^{nd}$ Law.

One very important point to note is the logical status of the 1st Law in physics: Its defense rests on (a) assuming that it is true and considering what this assumption implies; and (b) failing to find counterexamples to these implications in experiments. Just so, here we *assume* the applicability of an analogous Conservation Principle in the form suggested by Fermi, and attempt to articulate it as rigorously as possible. Ultimately the rationale for assuming that this principle applies to evidential systems will inhere in its utility.

**(1)** *Preliminaries* We have a data set $D_1$ and the graph of the LR corresponding to $D_1$ (and the hypotheses of interest), which we'll denote LR(A). LR(A) represents the "system" in its initial state. We are concerned with the effects on this system of a second data set, $D_2$. Let's designate the final state of the system, after consideration of $D_2$, as LR(B). Note that at this point we have not said anything specific about how we get from LR(A) to LR(B). For expository reasons, let's assume we're working with binomial data (n,x) and the (one-sided) binomial likelihood ratio as described in the main text.

We want to characterize the change in the LR graph from LR(A) to LR(B), or $\Delta$LR, in terms of change in evidence. Prior to defining evidence formally, we note that evidence measures a property of some key feature of the system. Let's denote this feature by $U_E$. The single most important feature of $U_E$ to note for the moment is that it must be a state variable associated with the LR graph. That is, for any given set of data (n, x), $U_E$ depends only on the data used to draw the LR graph, and not, for example, on anything related to the history of data collection. Otherwise, the evidence would depend upon things other than the LR graph, which by stipulation is precluded (see p. 6 above). As with its physical counterpart internal energy U, $U_E$ is defined only up to an additive constant, determined by arbitrarily selecting some standard state $U_E(0)$ against which other states are measured. (See Fermi (1956, p. 12-14) which refers to the indeterminate constant as "an essential feature of the concept of energy," so that only differences in energy are meaningfully quantifiable. This property is also essential to the log likelihood, which is closely related to the expected Shannon information, further supporting our proposal that (physical) energy and (evidential) information are mathematical counterparts of one another.)

We now seek to formally characterize the change in $U_E$ from state A to state B, or $\Delta U_E$, corresponding to $\Delta LR$. We know that the binomial LR graph for given n, x can be uniquely specified in terms of two quantities, but there is leeway regarding which two we choose. Ordinarily, we simply use the underlying (n, x) pair itself. But here we are looking for a more general theory, one that will allow us to measure E when the data do not take binomial form, so we want to work in terms of characteristics of the graph itself. Still, we have our choice of pairs of quantities we can use for this, e.g.: the area under the LR curve (let's call this $V_E$) and the maximum height; $V_E$ and the maximizing value of x/n; the maximizing x/n and curvature of the graph around its maximum; etc.. Any two of these uniquely specifies the corresponding graph for given n and x. Since at some level it does not matter which two quantities we pick, we elect to use $V_E$ and an as yet unspecified second quantity, which we'll call $P_E$. The only requirement for $P_E$ at this point is that $(V_E, P_E)$ together need to uniquely determine an admissible state of the LR graph, that is, the binomial LR for given (n, x). (Other forms of the likelihood may require different numbers of variables for unique specification of the LR graph. The circumstances under which this is the case, and the additional quantities required in these cases, will be a matter of practical importance and a subject of future investigation. But nothing fundamental to the concepts developed in the remainder of this document hinges on this point.)

**(2)** *Paths* We illustrate the notion of evidential paths (or simply "paths") with a concrete example. We begin with $D_1 = (n, x) = (4, 0)$ and $D_2 = (2, 1)$. Thus LR(A) corresponds to $(4, 0)$ while LR(B) corresponds to $(6, 1)$. We now consider the following two possibilities: the 5th toss $= 1$ and the 6th toss $= 0$; or the 5th toss $= 0$ and the 6th toss $= 1$. That is, the transformation from LR(A) to LR(B) could go through the intermediate states $(5,1)$ or $(5,0)$. We refer to these as Path1 and Path2 respectively. Figure 4 shows the corresponding LR graphs. As above, we can fully specify each LR graph in terms of our quantities $V_E$ and $P_E$. Thus, using a particular choice of $P_E$ (here we apply Eq. 7, Table 1 in main text) we can also represent each of the two paths on a $P_E V_E$ diagram, as shown in Figure 5. (A different choice of $P_E$ would change the numerical values of quantities computed below, but not the underlying principles.)

**(3)** *Work* We now need a word for the effect of $D_2$ on the LR graph in terms of V and P. Let's call this effect "work" ($-W_E$). (As in physics, the sign on $W_E$ is conventional, and here reflects the decision to describe the transformation produced by the new data as work done by the new data on the LR graph.) As Figure 5 illustrates, we need a measure of work that preserves the path-dependent nature of the transformation. In principle, we could consider either $\Delta W_E = V_E dP_E$ or $\Delta W_E = P_E dV_E$; we choose to work with the latter, for reasons having to do with our choice of $V_E$ and its relationship to physical volume, e.g., in terms of thermoscopic properties. Thus for our transformation from state A to state B, we have $-W_E = \int_A^B P_E dV_E$.

To illustrate, we can calculate $-W_E$ for each of our paths (from the numerical example above) as the sum of the work done at each step on each of the two paths, respectively. This gives us $-W_E(\text{Path1}) = 1.79$ and $-W_E(\text{Path2}) = 1.47$. Thus the work associated with the transformation from (4, 0) to (6, 1) depends on the order of the 2 intervening tosses. Or as Fermi (1956, p. 15) puts it, "the work performed [in going from the initial to final states] depends on whether we go by means of the first way or by means of the second way."

**(4) *Heat*** In the illustration above, since both paths begin and end with identical LR graphs, our state variable $U_E$ must change by the same amount in both cases. Our conservation principle instructs us that the change in the LR graph should reflect all and only the effects of $D_2$, but if we rest with the notion that $-W_E$ captures *all* of the effects of $D_2$, then we are left with the contradictory conclusion that $D_2$'s effects are both path-dependent and path-independent. There is only one way to reconcile these facts while adhering to the conservation principle, and that is to postulate the existence of an additional quantity. Let's call this additional quantity $Q_E$. While Q could simply stand for "quantity," because in thermodynamics it is called "heat," we refer to $Q_E$ as heat as well.

As long as we consider (physical) work as a fundamentally mechanical concept, it may be difficult to imagine a rationale for assuming an evidential, or purely mathematical, counterpart. In physics, it is convenient to illustrate the $1^{st}$ Law by contrasting mechanical work with non-mechanical forms of energy transfer, undoubtedly because mechanical work is so easily grasped and measured, while other concepts of energy, including heat itself, are strikingly abstract. However, from a formal point of view, what is essential to the Conservation Principle is not the equivalence between non-mechanical and mechanical forms of energy, but rather, the fact that in the face of path-dependent energy transfer, we need to have (at least) two forms of energy available to us to fully describe the desired dynamics.

**(5) *The First Law*** We can now formally articulate the conservation principle with which we started in precise notation, elevating it to the status of an evidential $1^{st}$ Law: $\Delta U_E = -W_E + Q_E$, or equivalently, $Q_E = \Delta U_E + W_E$. Note that we are able to express our Conservation Principle in this form in part because we've stipulated that $Q_E$ is the sole "something else" upon which the relationship between $U_E$ and $W_E$ depends. Whether or not this assumption is warranted will be put to the test by experimentation down the road.

One noteworthy feature of assuming the Evidential $1^{st}$ Law is that it permits us to calculate the quantities $\Delta U_E$ and $Q_E$, despite the utter abstractness of their definitions. Doing so, however, requires one further assumption, viz., specification of the mathematical relationships among $P_E$, $V_E$ and evidence E. To illustrate we assume that this relationship adheres to the ideal gas EqS as in the main text (see also Appendix B). We can then calculate $dU_E = C_V dE = -1.365$, which must be the same for both paths since E is a state variable. Plugging this number into the evidential 1st Law we obtain $Q_E$(Path1) = -3.1550 and $Q_E$(Path2) = -2.8350.

**(6) *The $2^{nd}$ Law***  In the main text we also defined a quantity $S_E$ (evidential entropy), and assumed $dS_E = \delta Q_E / E$. The $2^{nd}$ Law of Thermodynamics stipulates that in an isolated system, dS > 0 for "irreversible processes" and dS = 0 for reversible processes. It is readily confirmed that in the absence of energy influx, our evidential dS = 0. Recall too that all evidential transformations defined thus far are inherently reversible. It is unclear based on the above what the evidential analogue of irreversible processes would be, if any, and for this reason we postpone further formal development of an evidential 2nd Law. Nevertheless, we point out that some version of the 2nd Law is inevitable if the current formalism is to maintain coherence. In particular, it might appear unclear on the face of it why evidential energy could not "flow" equally well from a system at higher E to one at lower E or vice versa. But the same reasoning that precludes this in thermodynamics applies here: if information could flow equally in either direction, we could build an evidential perpetual motion machine by linking evidential Carnot engines (Fermi 1956 (orig. 1936)). Thus it must be the case that evidential energy "flows" only in one direction, and equivalently, that the only spontaneous changes in $S_E$ that can occur are increases. See (Vieland 2013) for further discussion.

**APPENDIX B: Relationship between E and Observed Fisher Information**  We note first that when $x/n$ is $<< 0.5$, the LR in the region $x/n > 0.5$ becomes small, so that the area under the LR in that region is negligible. Therefore, we can approximate $V_E$ (Eq. 3 in Table 1) by integrating from 0 to 1 (rather than from 0 to 0.5), which simplifies the arithmetic. In this case we have

$$V_E \approx \int_0^1 LR(\theta; n, x)\, d\theta = \int_0^1 2^n \theta^x (1-\theta)^{n-x}\, d\theta$$
$$= 2^n B(x+1, n-x+1) = 2^n \frac{x!(n-x)!}{(n+1)!} \tag{A1}$$

where $B$ is the beta function. Additionally, we apply a form of Stirling's approximation (Zwillinger),

$$\log(n!) \approx n \log n - n + \log \sqrt{2\pi n}. \tag{A2}$$

This gives

$$V_E \approx 2^n \frac{x!(n-x)!}{(n+1)!}$$
$$\approx 2^n \frac{x^x (n-x)^{n-x}}{(n+1)^{n+1}} \frac{e^{n+1}}{e^x e^{n-x}} \sqrt{\frac{2\pi x(n-x)}{n+1}}$$
$$= 2^n \frac{x^x (n-x)^{n-x}}{(n+1)^{n+1}} e \sqrt{\frac{2\pi x(n-x)}{n+1}}. \tag{A3}$$

Setting $C_V = \frac{1}{2}R$ (rather than $(3/2)R$ as used in the main text), and substituting the final expression in A3 into the expression for E (Eq. 6, Table 1), we obtain

$$E = \left(\frac{e^{S_E}}{V}\right)^2 \approx \left(\frac{2^n \left(\dfrac{x}{n}\right)^x \left(1-\dfrac{x}{n}\right)^{n-x}}{2^n \dfrac{x^x(n-x)^{n-x}}{(n+1)^{n+1}} e \sqrt{\dfrac{2\pi x(n-x)}{n+1}}}\right)^2$$

$$= \left(\left(\frac{n+1}{n}\right)^{2n} \frac{(n+1)^3}{e^2 \, 2\pi x(n-x)}\right). \tag{A4}$$

Noting that

$$\lim_{n \to \infty} \left(\frac{n+1}{n}\right)^n = e \tag{A5}$$

we see that as $n$ goes to infinity,

$$E \to \frac{(n+1)^3}{2\pi x(n-x)} = E_{APPROX}. \tag{A6}$$

For the Binomial distribution, the Observed Fisher Information, evaluated at the maximizing value $\theta = x/n$, is

$$FI_{OBS} = \frac{n^3}{x(n-x)}. \tag{A7}$$

Thus in large samples, we have

$$E \approx \frac{1}{2\pi} \times FI_{OBS}, \tag{A8}$$

where the approximation improves as $x/n$ approaches 0 and $n$ goes to infinity. Figure 6 illustrates the relationship between E (Eq. 6 in Table 1), $E_{APPROX}$ (from A6), and $FI_{OBS}/(2\pi)$ (from A7 and A8) for the model developed in the main text. Note that for fixed $x/n$, $FI_{OBS}$ is simply a constant times $n$. Thus E is approximately linear in $n$ for small fixed values of $x/n$, or equivalently, at low pressure $P_E$. This accords with the common use of the simple ln LR as an evidence measure in the statistical literature, a quantity which is itself linear in $n$ for fixed $x/n$. (However, E does not behave like the ln LR at higher $P_E$, but rather adapts to the full dynamics of the system described in the main text.) This suggests the possibility of calibrating units of E against the ln LR in low fixed-pressure settings, much the way the size of the degree Kelvin was set by calibrating against the air thermometer. It is also noteworthy that E converges to $FI_{OBS}/(2\pi)$ under exactly those circumstances in which real gases behave like ideal gases: large $n$, or equivalently, high E (the counterpart of high temperature in physical systems), and small $x/n$, or equivalently, low $P_E$ (the counterpart of low pressure in physical systems).

The approximation inherent in A1 itself improves as $x/n$ becomes smaller and $n$ becomes larger. Indeed, had we started from the outset with a two-sided hypothesis contrast (coin is fair vs. coin is not fair, without stipulating the direction of bias), the final expression in A1 would be exact rather than approximate, regardless of $n$. In this connection, what is most intriguing about Figure 6 is perhaps not where E and $FI_{OBS}$ coincide, but where they diverge: the relationship decays in the vicinity of the TrP(E), beyond which (reading left to right) E increases while $FI_{OBS}$ continues to decrease. Another way to express this is to say that E is taking explicit account of the particular hypothesis contrast (one-sided vs. two-sided) in a way that FI does not. We believe that this behavior can be related to phase transitions, following perhaps work such as (Prokopenko, Lizier et al.), which explicitly expresses statistical mechanical representations of thermodynamic systems in information theoretic terms. The connection between E and the observed Fisher Information - a relationship that we discovered only *after* development of E - also suggests a relationship between this work and (Frank), which posits the centrality of the Fisher Information itself to the dynamics of biological systems.

**ACKNOWLEDGMENTS** We thank the reviewers for their careful reading and insightful critical comments, which led to substantial improvements in the manuscript.

**Table 1. Fundamental Equations**

| Eq. # | Equation | Description |
|---|---|---|
| Eq. 1 | $LR = LR(\theta \mid n, x) = \dfrac{\theta^x (1-\theta)^{n-x}}{\frac{1}{2}^n}$ | Binomial likelihood ratio, for $n$ tosses, $x$ heads, P(heads)=$\theta$, where the denominator hypothesis is $\theta = \frac{1}{2}$ (no bias) |
| Eq. 2 | $PV = RT$ | Ideal gas equation of state (EqS) for 1 mole, where P=pressure, V=volume, R is a positive (universal scaling) constant, and T= absolute temperature |
| Eq. 3 | $V_E = \int LR(\theta) \, d\theta = \int 2^n \theta^x (1-\theta)^{n-x} \, d\theta$ | Evidential Volume |
| Eq. 4 [*] | $S = C_V \log T + R \log V$ | Alternative form of EqS together with $1^{st}$ and $2^{nd}$ Laws of Thermodynamics, where S=entropy and $C_V$ is a positive constant. |
| Eq. 5 | $S_E = n\left[\dfrac{x}{n}\log\left(\dfrac{x}{n}\right) + \dfrac{n-x}{n}\log\left(1-\dfrac{x}{n}\right) + \log 2\right] + k$ | Evidential Entropy, defined up to a constant $k$ |
| Eq. 6 | $E = \dfrac{\exp\left\{\frac{S_E}{C_V}\right\}}{V_E^{R/C_V}}$ | Evidence |
| Eq. 7 | $P_E = \dfrac{RE}{V_E} = R\left(\dfrac{\exp\left\{\frac{S_E}{C_V}\right\}}{V_E^{(R/C_V + 1)}}\right)$ | Evidential Pressure |
| Eq. 8 | $W_E = \int\limits_A^B P_E \, dV_E$ | Evidential Work |

[*] Note that Eq. 4, while a standard representation in thermodynamics texts, is not adequate for expressing the behavior of ideal gasses at very low temperatures, as T approaches its true minimum value of 0 [see Fermi (1956), p. 147]. Hence we avoid this part of the range throughout this paper.

**Figure Legends**

Fig. 1 Behavior of binomial evidential system: (a) E as a function of *x/n* for various values of *n*; (b) *x/n* at TrP(E) as a function of *n*; (c) E as a function of (*n, x*).

Fig. 2 (a) Various "isothermal" contours of the binomial evidential system as a function of (*n, x*); (b) These same isotherms plotted as a function of (*n, x/n*); (c) The E=2.25 isotherm on a $P_E V_E$ diagram. Starting from the TrP: Moving to the left in panels (a) and (b) corresponds to moving to the right in panel (c), and vice-versa, moving to the right in (a) and (b) corresponds to moving to the left in (c).

Fig. 3 Two evidential Carnot cycles for the binomial system in the format of a classic PV diagram. Cycle 1 (on left) operates between $E_1$=1 and $E_2$=2; cycle 2 (on right) operates between $E_1$=2 and $E_2$=4. Both cycles start in the upper left and travel clockwise. For each stroke *i* (*i*=A,B,C,D; see text), work $W_i$ is calculated as the area beneath the line connecting the nodes. For each cycle, the ratio $W_C/W_A$ equals the ratio of evidence levels, or ½, as indicated visually by the right-hand panel for each cycle. This illustrates the general principle that for any cyclic transformation taking a system from evidence E to (½)E, the ratio of mechanical work performed at the two evidence levels is the same.

Fig. 4 Four LR graphs corresponding to the initial and final states of the system, as well as the two intermediate states of the system corresponding to Path1 and Path2, respectively. Each path is labeled with its corresponding (n, x) pair.

Fig. 5 $P_E V_E$ plots corresponding to (a) Path1 and (b) Path2. The labels show the corresponding (n, x) pairs for the initial, intermediate and final states of the system, respectively.

Fig. 6 E, $E_{APPROX}$, and $FI_{OBS}/(2\pi)$ for n = 10, 20, 50, and 100 (panels (a)-(d) respectively). The three quantities converge towards one another as *n* increases to the left of the TrP(E). Note that at *x/n* = 0, E has a finite value while both $E_{APPROX}$ and $FI_{OBS}/(2\pi)$ become infinite.

Figure 1

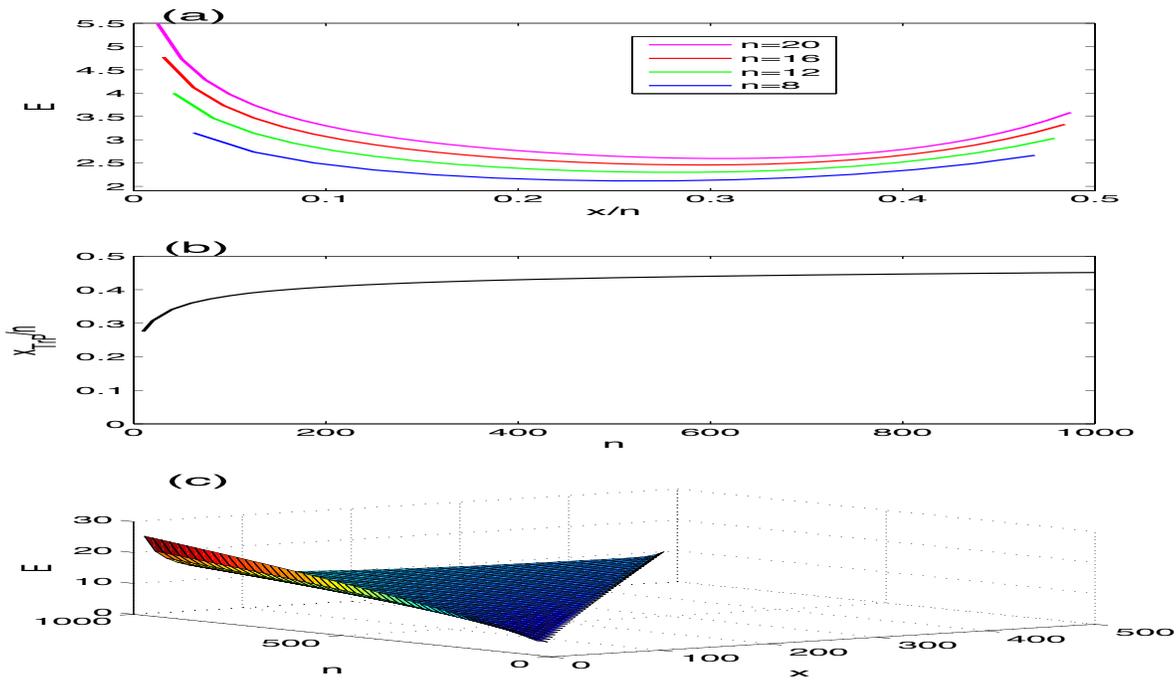

Figure 2

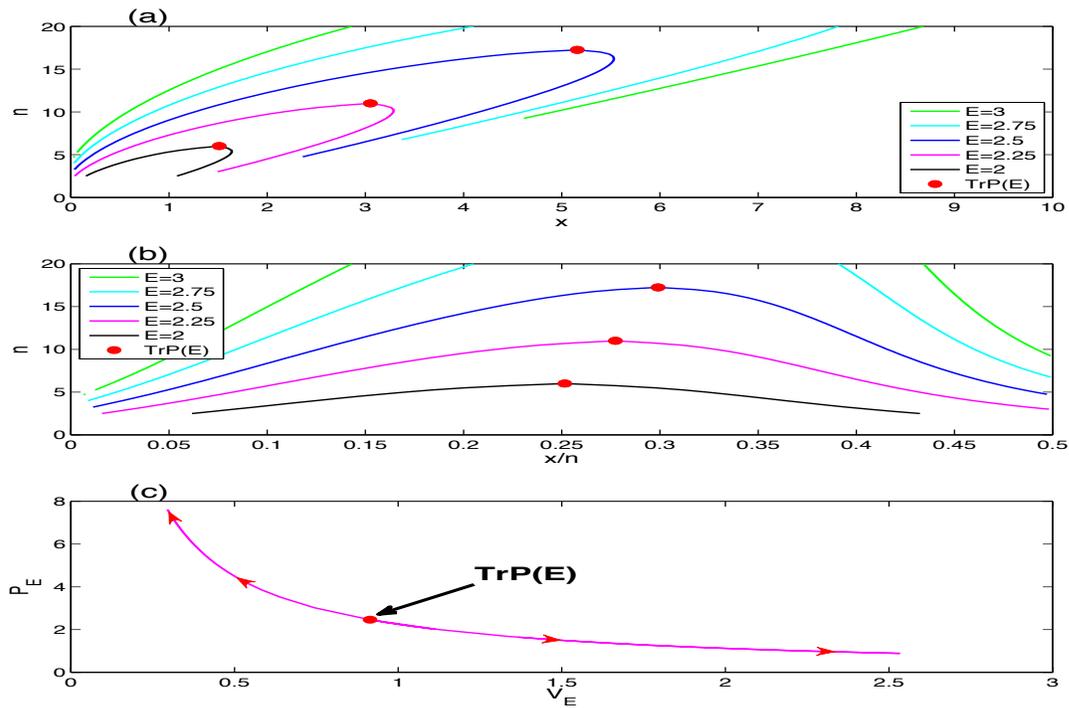

Figure 3

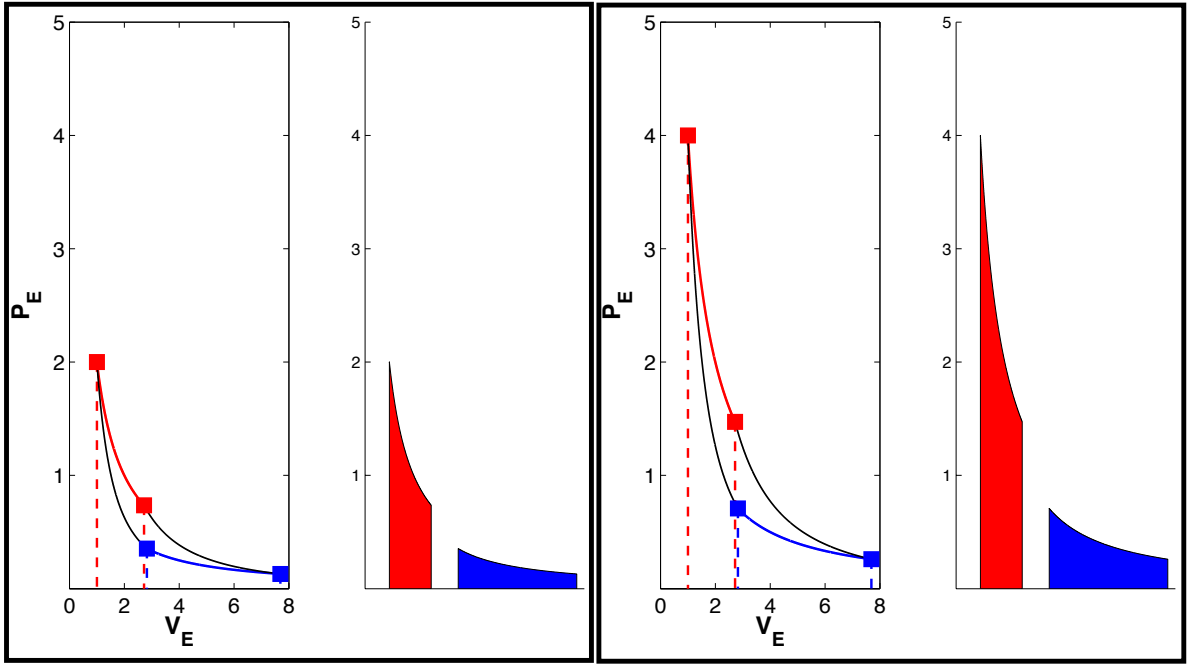

Figure 4

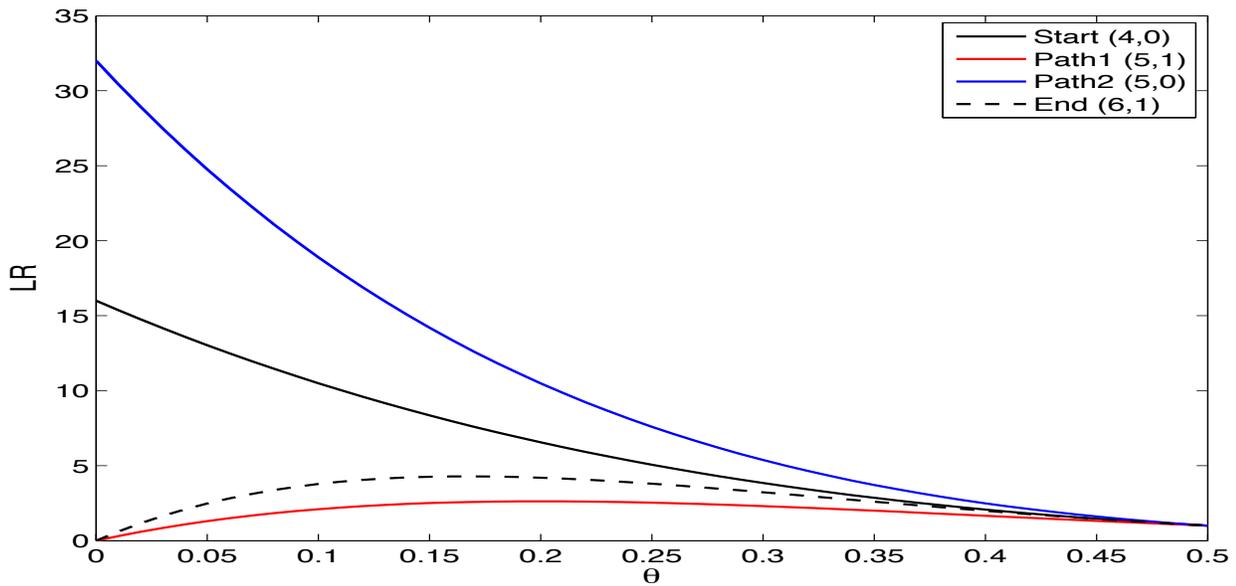

Figure 5

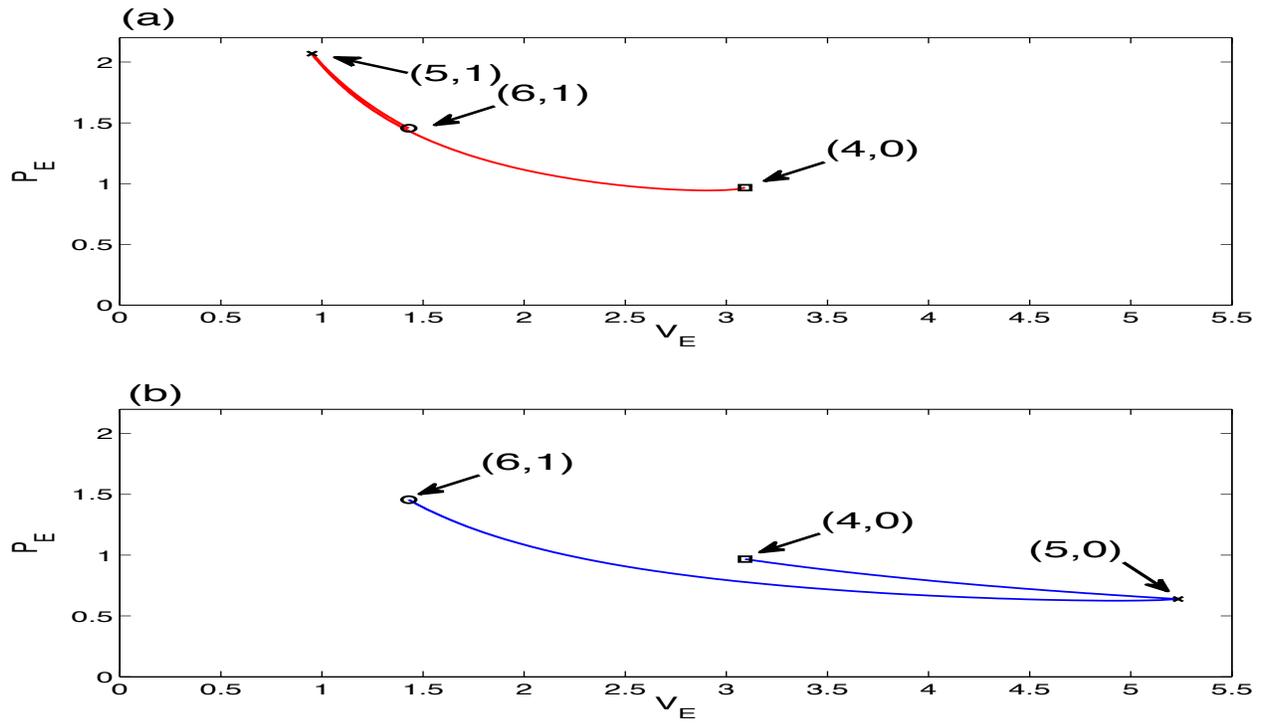

Figure 6

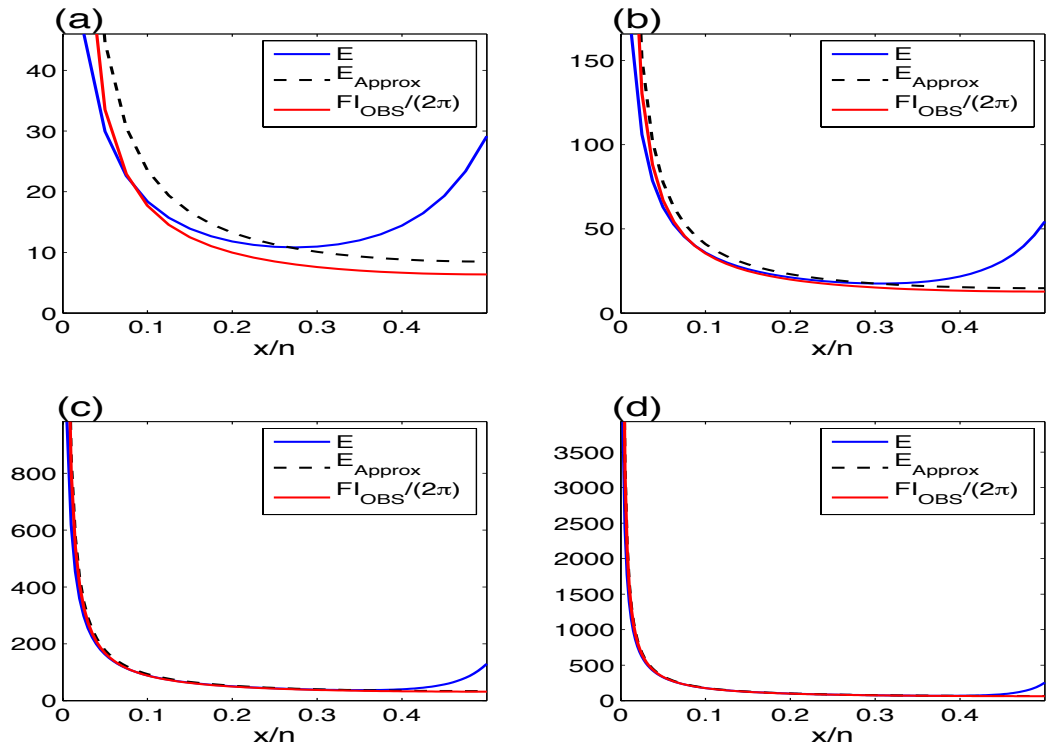